\documentclass[11pt,reqno]{amsart} % please use amsart at 11pt
\usepackage{amssymb,latexsym}
\usepackage[draft=true]{hyperref}
\usepackage{cite} % to get Refs. [1,2,3] typeset as [1--3] automatically
\usepackage{amssymb}
\usepackage{mathrsfs}
\usepackage{amsmath}
\usepackage{amsfonts}
\usepackage{pifont}
\usepackage{graphicx,amssymb,mathrsfs,amsmath}
\usepackage{tocvsec2}
\usepackage{color}
\usepackage{graphicx}
\usepackage{float}
\usepackage{epstopdf}
\usepackage{amssymb, amsmath, amsthm, graphics, comment, xspace, enumerate, xcolor}
\usepackage[height=190mm,width=130mm]{geometry} % this is the journal
\theoremstyle{plain}
\newtheorem{theorem}{Theorem}
\newtheorem{lemma}{Lemma}

\theoremstyle{definition}

\theoremstyle{remark}
\newtheorem{remark}{Remark}

\numberwithin{equation}{section} % to get equations numbered
\begin{document}

\title[New refinements of Cusa-Huygens inequality]{New refinements of Cusa-Huygens inequality} % please provide
                                % an abbreviated title

\author[C. Chesneau]{Christophe Chesneau}
\address{LMNO\\ University of Caen-Normandie\\ Caen\\ France}
\email{christophe.chesneau@unicaen.fr}

\author[M. Kosti\'{c}]{Marko Kosti\' c}
\address{Faculty of Technical Sciences\\
University of Novi Sad\\
Trg D. Obradovi\' ca 6\\ 21125 Novi Sad\\ Serbia}
\email{marco.s@verat.net}

\author[B. Male\v{s}evi\'{c} ]{Branko Male\v{s}evi\'{c}}
\address{School of Electrical Engineering\\
University of Belgrade, Serbia\\}
\email{branko.malesevic@etf.bg.ac.rs}

\author[B. Banjac]{Bojan Banjac}
\address{Faculty of Technical Sciences\\
University of Novi Sad, Serbia\\}
\email{bojan.banjac@uns.ac.rs}

\author[Y. J. Bagul]{Yogesh J. Bagul}
\address{Department of Mathematics\\ K. K. M. College Manwath\\ Dist: Parbhani(M.S.) - 431505\\
India}
\email{yjbagul@gmail.com}

%\thanks{The third author was supported in part by Ministry
%of Science and Technological Development, Republic of Serbia, Grant No.~174024.}

\begin{abstract}
In the paper, we refine and extend Cusa-Huygens inequality by simple functions. In particular, we determine sharp bounds for $\sin(x) /x$ of the form $(2+\cos(x))/3 -(2/3-2/\pi)\Upsilon(x)$, where $\Upsilon(x)>0$ for $x\in (0, \pi/2)$, $\Upsilon(0)=0$ and $\Upsilon(\pi/2)=1$, such that $\sin x/x$ and the proposed bounds coincide at $x=0$ and $x=\pi/2$. The hierarchy of the obtained bounds is discussed, along with a graphical study. Also, alternative proofs of the main result are given. 
\end{abstract}

%\dedicatory{This paper is dedicated to Professor X on his 125th birthday.}



\maketitle

\section{Introduction}\label{intro}
The Cusa-Huygens inequality \cite{huygens, mitrinovic} is one of the celebrated inequalities in the theory of analytic inequalities involving trigonometric  functions. It is stated as follows:
\begin{align}\label{eqn1.1}
\frac{\sin x}{x}<\frac{2+\cos x}{3}, \quad x\in (0,\pi/2).
\end{align}
In the recent paper \cite[Remark 4.1]{zhu}, it is remarked that the inequality (\ref{eqn1.1}) is, in fact, true for all $ x \neq 0. $ Due to the symmetry of the functions involved at both sides of (\ref{eqn1.1}), it suffices to consider the inequality on the right half of the real line. Although (\ref{eqn1.1}) is true for $ x > 0 $ and it is not sharp everywhere. So the interest among researchers has been to consider the inequality in the natural intervals $ (0, \pi/2) $ or $ (0, \pi). $ In \cite{zhu}, L. Zhu obtained Cusa-Huygens type inequalities on a wider range $ (0, \pi). $ We must emphasize here that the sharp Cusa-Huygens type inequalities on  a wider range $ (0, \pi) $  have also appeared in \cite{bagul2, dhaigude, malesevic1, mortici}.

In \cite{chen, sandor}, the following Cusa-Huygens type double inequality  was established:
\begin{align}\label{eqn1.2}
\left(\frac{2+\cos x}{3}\right)^{\alpha}<\frac{\sin x}{x}<\left(\frac{2+\cos x}{3}\right)^{\zeta}, \quad x\in (0,\pi/2),
\end{align}
where $\alpha =\ln (\pi/2)/\ln(3/2) \approx 1.11374$ and $\zeta=1$ are the best possible constants. Simple alternative  proofs of \eqref{eqn1.2} are offered in \cite{bagul, bagul1}. For other details about inequalities \eqref{eqn1.1} and \eqref{eqn1.2}, we refer readers to \cite{bagul, bagul1, bagul2, bagul3, bagul4, chen, chen1, dhaigude, huygens, malesevic, malesevic1, mitrinovic, mortici, neuman, sandor, sandor1, zhu}. 

The main objective of this paper is to provide a manageable and sharp alternative to the double inequality (\ref{eqn1.2}) in the following sense: (i) in the applications, the exponentiated version of $(2+\cos(x))/3$ is not always manageable, so we have in mind to keep the overall simplicity of the former Cusa-Huygens inequality by taking the unit exponent, and (ii) the values of the bounds in (\ref{eqn1.2}) differ to the one of $\sin x /x$ when $x=\pi/2$, making them not sharp for an interval of values of $x$ close to $\pi/2$.  For these reasons, we focus our attention on sharp bounds of $\sin x/x$ of the following form: 
\begin{align}\label{r(x)}
Z(x):=\frac{2 + \cos x}{3} - \left(\frac{2}{3} - \frac{2}{\pi}\right)\Upsilon(x),
\end{align}
 where $\Upsilon(x)>0$ for $x\in (0, \pi/2)$, $\Upsilon(0)=0$ and $\Upsilon(\pi/2)=1$. Hence, a simple form for $\Upsilon(x)$ implies a tractable expression for $Z(x)$. Candidates of such functions are proposed and two theorems are proved. The importance of the finding is illustrated by a graphical study, showing the high degree of sharpness of the results. Also, the hierarchy between the obtained bounds are examined. Alternative proofs are also provided. 
 
The rest of the paper is divided into the following sections.  Section \ref{sec2} contains the main results. Section \ref{sec3} is devoted to the proofs of these results. Diverse complements are offered in Section \ref{sec4}, including a comparison of the proposed bounds, a graphical study and alternative proofs.

 \section{Main results}\label{sec2}
The following theorem  presents sharp bounds for $\sin x/x$ of the form \eqref{r(x)}, with $\Upsilon(x)$ defined as exponentiated version of $(\pi/2-1)^{-1}(x-\sin x)$. 
\begin{theorem}\label{thm1}
Let $ x \in (0, \pi/2). $ Then the double inequality
\begin{align}\label{eqn1.3}
\frac{2 + \cos x}{3} - \left(\frac{2}{3} - \frac{2}{\pi}\right)\Phi_1(x) < \frac{\sin x}{x} < \frac{2 + \cos x}{3} - \left(\frac{2}{3} - \frac{2}{\pi}\right)\Phi_2(x)
\end{align}
where $ \Phi_1(x) := (\pi/2-1)^{-1}(x-\sin x) $ and $ \Phi_2(x) := (\pi/2-1)^{-2}(x-\sin x)^2 $. 
\end{theorem}

In a similar fashion to Theorem \ref{thm1}, Theorem \ref{thm2} presents sharp bounds for $\sin x/x$ of the form \eqref{r(x)}, with $\Upsilon(x)$ defined as exponentiated version of $\sin x - x \cos x $.  
\begin{theorem}\label{thm2}
Let $ x \in (0, \pi/2). $ Then, the double inequality
\begin{align}\label{eqn1.4}
\frac{2 + \cos x}{3} - \left(\frac{2}{3} - \frac{2}{\pi}\right)\Psi_1(x) < \frac{\sin x}{x} < \frac{2 + \cos x}{3} - \left(\frac{2}{3} - \frac{2}{\pi}\right)\Psi_2(x)
\end{align}
where $ \Psi_1(x) := \sin x - x \cos x $ and $ \Psi_2(x) := (\sin x - x \cos x)^2 $. 
\end{theorem}

Since $\Psi_2(x) $ and $\Psi_2(x) $ are positive, the Cusa-Huygens inequality in (\ref{eqn1.1}) is clearly refined. Moreover, as linear combinations of simple functions, the obtained bounds are quite manageable for diverse analytical manipulations involving integration, series, and so on. 

The main proofs of these theorems are based on several general and specific results described in the next section. Some of them can be of independent interest, for aims beyond the scope of this study. Let us mention that the sharpness of the obtained bounds, as well as alternative proofs, will be discussed later. 
\section{Proofs}\label{sec3}
\subsection{Preliminaries and lemmas}

The following series expansions can be found in \cite[1.411(7, 11)]{grad}:
\begin{align}\label{eqn2.1}
\cot x = \frac{1}{x} - \sum_{k=1}^{\infty} \frac{2^{2k}\vert B_{2k} \vert}{(2k)!}x^{2k-1}, \quad  x \in (0, \pi)
\end{align}
and
\begin{align}\label{eqn2.2}
\frac{1}{\sin x} = \frac{1}{x} + \sum_{k=1}^{\infty}\frac{2(2^{2k-1}-1)\vert B_{2k} \vert}{(2k)!}x^{2k-1}, \quad  x \in (0, \pi).
\end{align}
From (\ref{eqn2.2}), we can write
\begin{align}\label{eqn2.3}
\frac{x}{\sin x} = 1 + \sum_{k=1}^{\infty}\frac{2(2^{2k-1}-1)\vert B_{2k} \vert}{(2k)!}x^{2k}, \quad  x \in (0, \pi).
\end{align}
Similarly, we obtain
\begin{align}\label{eqn2.4}
\left(\frac{x}{\sin x}\right)^2 = -x^2 (\cot x)' = 1 + \sum_{k=1}^{\infty}\frac{2^{2k}(2k-1) \vert B_{2k} \vert}{(2k)!}x^{2k}, \quad  x \in (0, \pi)
\end{align}
and
\begin{align}\label{eqn2.5}
\frac{\cos x}{\sin^2 x} = \left(\frac{1}{\sin x}\right)' =  \frac{1}{x^2} - \sum_{k=1}^{\infty}\frac{(2k-1)(2^{2k}-2) \vert B_{2k} \vert}{(2k)!}x^{2k-2}, \quad x \in (0, \pi),
\end{align}
respectively, from (\ref{eqn2.1}) and (\ref{eqn2.2}).

For Lemma \ref{lem1}, we refer to \cite{anderson}. Lemma \ref{lem1} is known as l'H\^{o}pital's rule of monotonicity.

\begin{lemma}\label{lem1}
 Let $ f, g : [m, n] \rightarrow \mathbb{R} $ be two continuous functions which are differentiable on $(m, n)$ and $ g'(x) \neq 0 $ in $ (m, n)$. If $ f'(x)/g'(x) $ is increasing (or decreasing) on $ (m, n),$ then the functions $(f(x) - f(m))/(g(x) - g(m))$ and $(f(x) - f(n))/(g(x) - g(n))$ are also increasing (or decreasing) on $ (m, n). $
If $ f'(x)/g'(x) $ is strictly monotone, then the monotonicity in the conclusion is also strict.
\end{lemma}

\begin{lemma}(\cite[Lemma 2]{bagul3})\label{lem2} The function
$$
P(x):=\frac{x^2 \sin x}{\sin x-x \cos x},\quad x\in (0,\pi/2)
$$
is positive and strictly decreasing on $(0,\pi/2)$.
\end{lemma}

Lemma \ref{lem3} below can be found in \cite{heikkala, alzer}.

\begin{lemma}\label{lem3}
Let $ A(x) := \sum_{n = 0}^{\infty}a_n x^n $ and $ B(x) := \sum_{n = 0}^{\infty}b_n x^n $ be convergent for $ \vert x \vert < T ,$ where $ a_n $ and $ b_n $ are real numbers for $ n = 0, 1, 2, \cdots $ such that $ b_n > 0$, and $T>0$ is a fixed constant.  If the sequence ${a_n/b_n}$ is strictly increasing(or decreasing), then the function $ A(x)/B(x) $ is also strictly increasing(or decreasing) on $ (0, T).$
\end{lemma}

The next lemma is about lower and upper bounds for a ratio involving absolute Bernoulli numbers. It is established in \cite{qi}.
\begin{lemma}\label{lem4}
For $ k \in \mathbb{N}, $ the Bernoulli numbers satisfy
$$ \frac{(2^{2k - 1} - 1)}{(2^{2k + 1} - 1)} \frac{(2k + 1) (2k + 2)}{\pi^2} < \frac{\vert B_{2k + 2} \vert}{\vert B_{2k} \vert} <  \frac{(2^{2k} - 1)}{(2^{2k + 2} - 1)} \frac{(2k + 1) (2k + 2)}{\pi^2}. $$
\end{lemma}
The following result is for future technical considerations. 
\begin{lemma}\label{lem5}
The function $$ Q(x) := \frac{3 \sin x - 2x - x \cos x}{x^5}, \quad x \in (0, \pi/2) $$ is negative and strictly increasing on $ (0, \pi/2).$
\end{lemma}
\begin{proof}
First, $ Q(x) $ is negative due to inequality (\ref{eqn1.1}).
Now, let $$ q_1(x) := 3 \sin x - 2x - x \cos x, \quad  h_1(x) := x^5, $$
$$ q_2(x) := x \sin x + 2\cos x -2,  \quad h_2(x) := 5 x^4$$ and
$$ q_3(x) := x \cos x - \sin x, \ h_3(x) := 20 x^3.$$
Then $$ q_i(0+) = h_i(0+) = 0 \ (i = 1, 2, 3), \quad \frac{q_i'(x)}{h_i'(x)} = \frac{q_{i+1}(x)}{h_{i+1}(x)} \ (i = 1, 2) $$ and
$$ \frac{q_3'(x)}{h_3'(x)} = -\frac{1}{60} \frac{\sin x}{x}. $$ Since the function $ \sin x/x $ is strictly decreasing on $ (0, \pi/2), $ repeated application of Lemma \ref{lem1} gives that $ Q(x) $ is strictly increasing on $ (0, \pi/2).$
\end{proof}
The following result is also an important ingredient of a main proof. 
\begin{lemma}\label{lem6}
The function $$ T(x) := \frac{x^2 \sin x}{x - \sin x}, \quad x \in (0, \pi) $$ is positive and strictly decreasing on $ (0, \pi).$
\end{lemma}
\begin{proof}
Clearly, $ T(x) $ is positive since $ x > \sin x, \ x \in (0, \pi). $ Thus, it suffices to prove that $ 1/T(x) $ is strictly increasing on $ (0, \pi).$ We have
$$ \frac{1}{T(x)} = \frac{1}{x^2} \left(\frac{x}{\sin x} - 1 \right). $$ By (\ref{eqn2.3}), we get
$$ \frac{1}{T(x)} = \sum_{k=1}^{\infty}\frac{2(2^{2k-1}-1)\vert B_{2k} \vert}{(2k)!}x^{2k - 2}, \quad x \in (0, \pi)$$
which is obviously strictly increasing on $ (0, \pi).$
\end{proof}

\subsection{Proofs of main results} \ 

{\it Proof of Theorem \ref{thm1}.} We first prove the left inequality of \eqref{eqn1.3}. Let us set
$$ f(x) := \frac{\frac{\sin x}{x} - \frac{2 + \cos x}{3} }{x - \sin x} = \frac{3 \sin x - x \cos x - 2x}{3 x^2 - 3x \sin x}, \quad x \in (0, \pi/2). $$
As intermediary functions, we consider $$ g_1(x) := 3 \sin x - x \cos x - 2x, \quad  h_1(x) := 3x^2 - 3x \sin x, $$$$ g_2(x) := x \sin x + 2\cos x - 2, \quad h_2(x) := 6x - 3\sin x - 3x \cos x, $$ and $$  g_3(x) := x\cos x - \sin x, \quad  h_3(x) := 3x \sin x - 6\cos x + 6.$$  Then
$$ g_i(0+) = h_i(0+) = 0 \ (i = 1, 2, 3), \quad  \frac{g_i'(x)}{h_i'(x)} = \frac{g_{i+1}(x)}{h_{i+1}(x)} \ (i = 1, 2) $$ and
$$ \frac{g_3'(x)}{h_3'(x)} = \frac{-x \sin x}{9 \sin x + 3x \sin x} = -\frac{1}{3} \left(\frac{x}{3 + u(x)} \right) $$ with $$ u(x) := \frac{x}{\tan x}. $$ The function $ u(x) $ is positive and strictly decreasing on $ (0, \pi/2). $ Applying Lemma \ref{lem1} repeatedly, we get that $ f(x) $ is strictly decreasing on $(0, \pi/2). $ Hence, for $x < \pi/2$, we have $$ f\left(\frac{\pi}{2}-\right) < f(x).$$ Since $$ f\left(\frac{\pi}{2}-\right) = \frac{\frac{2}{\pi} - \frac{2}{3}}{\frac{\pi}{2} - 1 }, $$ we get required left inequality of \eqref{eqn1.3}.

Now we prove the right inequality of (\ref{eqn1.3}). Let
$$ F(x) := \frac{\frac{\sin x}{x} - \frac{2 + \cos x}{3}}{(x - \sin x)^2} = \frac{3 \sin x - x \cos x - 2x}{3 x^3 - 6 x^2 \sin x + 3x \sin^2 x}.  $$ It can also be written as
$$ F(x) = \frac{3 \frac{1}{\sin x} - \frac{2}{x}\left(\frac{x}{\sin x}\right)^2 - x \frac{\cos x}{\sin^2 x}}{3x\left(\frac{x}{\sin x}\right)^2 - 6x \left(\frac{x}{\sin x}\right)+3x}. $$ Utilizing \eqref{eqn2.2}, \eqref{eqn2.3}, \eqref{eqn2.4} and \eqref{eqn2.5}, we can write
{\footnotesize
\begin{align*}
& F(x) = \\
& \left(\frac{\sum_{k=1}^{\infty}\frac{6(2^{2k-1}-1)}{(2k)!} \vert B_{2k} \vert x^{2k-1} - \sum_{k=1}^{\infty} \frac{2^{2k+1}(2k-1)}{(2k)!} \vert B_{2k} \vert x^{2k-1} + \sum_{k=1}^{\infty} \frac{(2k-1)(2^{2k}-1)}{(2k)!} \vert B_{2k} \vert x^{2k-1}}{\sum_{k=1}^{\infty} \frac{3 \cdot 2^{2k}(2k-1)}{(2k)!} \vert B_{2k} \vert x^{2k+1} - \sum_{k=1}^{\infty} \frac{3 \cdot 2^2  (2^{2k-1}-1)}{(2k)!} \vert B_{2k} \vert x^{2k+1}}\right) \\
&= \frac{\sum_{k=1}^{\infty}2\left[3(2^{2k-1}-1)-2^{2k}(2k-1)+(2k-1)(2^{2k-1}-1)\right]\frac{\vert B_{2k} \vert }{(2k)!} x^{2k-1}}{\sum_{k=1}^{\infty} 12\left[(2k-1)2^{2k-2}-(2^{2k-1}-1)\right]\frac{\vert B_{2k} \vert }{(2k)!}x^{2k+1}} \\
&= \frac{\sum_{k=2}^{\infty}2\left[3(2^{2k+1}-1)-2^{2k+2}(2k+1)+(2k+1)(2^{2k+1}-1)\right]\frac{\vert B_{2k+2} \vert }{(2k+2)!} x^{2k+1}}{\sum_{k=1}^{\infty} 12\left[(2k-1)2^{2k-2}-(2^{2k-1}-1)\right]\frac{\vert B_{2k} \vert }{(2k)!}x^{2k+1}} \\
&= \frac{\sum_{k=2}^{\infty}4\left[-2k \cdot 2^{2k} + 2 \cdot 2^{2k} -k -2\right] \frac{\vert B_{2k+2} \vert }{(2k+2)!} x^{2k+1}}{\sum_{k=2}^{\infty} 3 \left[ 2k \cdot 2^{2k} - 3 \cdot 2^{2k} + 4 \right] \frac{\vert B_{2k} \vert }{(2k)!}x^{2k + 1}} \\
&= \left(\frac{-1}{\sum_{k=2}^{\infty}\frac{180 \vert B_{2k} \vert}{(2k)!}(2k \cdot 2^{2k} - 3 \cdot 2^{2k} + 4 )x^{2k-2}} := M(x) \right) \\
&- \left(\frac{\sum_{k=2}^{\infty}  \frac{4 \vert B_{2k+2} \vert}{(2k+2)!}(2k \cdot 2^{2k} - 2 \cdot 2^{2k} + k +2)x^{2k+1}}{\sum_{k=2}^{\infty} \frac{3 \vert B_{2k} \vert}{(2k)!}(2k \cdot 2^{2k} - 3 \cdot 2^{2k} +4) x^{2k+1}} := \frac{A(x)}{B(x)}\right).
\end{align*}}
Now $$ 2k \cdot 2^{2k} - 2 \cdot 2^{2k} + k + 2 > 0, \ (k= 2, 3, 4, \cdots), $$ and $$ 2k \cdot 2^{2k} -3 \cdot 2^{2k} + 4 > 0, \ (k = 2, 3, 4, \cdots).$$
Therefore, $ M(x)$ is strictly increasing on $ (0, \pi/2). $ Next, let us notice that
$$ \frac{A(x)}{B(x)} = \frac{\sum_{k=2}^{\infty} a_k x^{2k+1}}{\sum_{k=2}^{\infty}b_k x^{2k+1}} $$ where
$$ a_k := \frac{4 \vert B_{2k+2} \vert}{(2k+2)!}(2k \cdot 2^{2k} - 2 \cdot 2^{2k} + k +2) > 0 $$ and
$$ b_k := \frac{3 \vert B_{2k} \vert}{(2k)!}(2k \cdot 2^{2k} - 3 \cdot 2^{2k} +4) > 0. $$
So, $$ \frac{a_k}{b_k} = \frac{2}{3} \frac{\vert B_{2k + 2} \vert }{\vert B_{2k} \vert}\frac{1}{(k + 1)(2k + 1)}\frac{(2k \cdot 2^{2k} - 2 \cdot 2 ^{2k} + k + 2)}{(2k \cdot 2^{2k} - 3 \cdot 2^{2k} + 4)} := c_k $$ and
$$ c_{k+1} = \frac{1}{6} \frac{\vert B_{2k+4} \vert}{\vert B_{2k+2} \vert} \frac{1}{(k+2)(2k+3)} \frac{(8k \cdot 2^{2k} + k + 3)}{(2k \cdot 2^{2k} - 2^{2k} + 1)}. $$ Hence
\begin{align*}
 \frac{c_{k+1}}{c_{k}} = \frac{1}{4} \frac{\vert B_{2k + 4} \vert}{\vert B_{2k + 2} \vert} &\frac{\vert B_{2k} \vert}{\vert B_{2k + 2} \vert} \frac{(k+1)(2k + 1)}{(k + 2)( 2k + 3)} \\
&\times \frac{(8k \cdot 2^{2k} + k + 3)(2k \cdot 2^{2k} - 3 \cdot 2^{2k} + 4)}{(2k \cdot 2^{2k} - 2^{2k} + 1)(2k \cdot 2^{2k} - 2 \cdot 2^{2k} + k + 2)}.
\end{align*}
 By Lemma \ref{lem4}, we have $$ \frac{\vert B_{2k+4} \vert}{\vert B_{2k +2} \vert} < \frac{2(4 \cdot 2^{2k}-1)(k+2)(2k+3)}{\pi^2 (16 \cdot 2^{2k}-1)} $$ and
$$ \frac{\vert B_{2k} \vert}{\vert B_{2k + 2} \vert} < \frac{\pi^2 (2 \cdot (2^{2k} - 1)}{(2^{2k}-2)(k+1)(2k+1)}.  $$ Then
\begin{align*}
 4 \frac{c_{k+1}}{c_k} < &\frac{2(4 \cdot 2^{2k} - 1)(2 \cdot 2^{2k}-1)}{(16 \cdot 2^{2k} - 1)(2^{2k}-2)} \\
 &\times \frac{(8k \cdot 2^{2k} + k + 3) (2k \cdot 2^{2k} - 3 \cdot 2^{2k} +4)}{(2k \cdot 2^{2k} - 2^{2k} +1)(2k \cdot 2^{2k} - 2 \cdot 2^{2k} + k +2)} := N(k).
\end{align*}
We claim that $$ N(k) < 4, $$ i.e.,
\begin{align*}
2(4 &\cdot 2^{2k} - 1)( 2 \cdot 2^{2k} - 1)(8k \cdot 2^{2k} + k + 3)(2k \cdot 2^{2k} -3 \cdot 2^{2k} + 4) \\
 &< 4(16 \cdot 2^{2k} -1)(2^{2k} - 2)(2k \cdot 2^{2k}-2^{2k} + 1)(2k \cdot 2^{2k} -2 \cdot 2^{2k} + k + 2),
\end{align*}
or 
\begin{align*}
&(8 \cdot 2^{4k} - 6 \cdot 2^{2k} +1) \\
&\times (16k^2 \cdot 2^{4k} - 24k \cdot 2^{4k} + 2k^2 \cdot 2^{2k} + 35k \cdot 2^{2k} - 9 \cdot 2^{2k} + 4k + 12) \\
&< (32 \cdot 2^{4k} - 66 \cdot 2^{2k} + 4) \\
&(4k^2 \cdot 2^{4k} - 6k \cdot 2^{4k} + 2 \cdot 2^{4k} + 2k^2 \cdot 2^{2k} + 5k \cdot 2^{2k} - 4 \cdot 2^{2k} + k + 2).
\end{align*}
Equivalently, after some computations, we get
\begin{align*}
(&120k^2 + 188)\cdot 2^{6k} + (120k^2 + 120)\cdot 2^{4k} + \left[(57k + 67)\cdot 2^{2k} + 4\right] \\
&< (64\cdot 2^{2k})\cdot 2^{6k} + (132k \cdot 2^{2k}) \cdot 2^{4k} + \left[(186\cdot 2^{2k} + 6k^2)\cdot 2^{2k} \right]
\end{align*}
which is clearly true for $ k = 2, 3, 4, \cdots. $ This implies that $ c_{k+1} < c_k ,$ i.e., a sequence $ \left\lbrace a_k/b_k\right\rbrace_{k=2}^{\infty} $ is strictly decreasing. By Lemma \ref{lem3}, $A(x)/B(x) $ is strictly decreasing on $ (0, \pi/2) $ giving us that $ F(x) $ is strictly increasing on $ (0, \pi/2). $ Consequently, for $x < \pi/2 $, 
$$ F(x) < F\left(\frac{\pi}{2}-\right)  = \left(\frac{2}{\pi} - \frac{2}{3}\right)\left(\frac{\pi}{2} - 1 \right)^{-2} $$ and
the desired right inequality of \eqref{eqn1.3} follows. \hfill  $ \qed $

\begin{remark}Fortunately, we can give a very simple proof of the right inequality of \eqref{eqn1.3}. Indeed, we  can write $ F(x) $ as follows:
$$ F(x) = \frac{1}{3} \cdot Q(x) \cdot R(x) \cdot \left[T(x)\right]^2,  $$
where $$ Q(x) := \frac{3 \sin x - 2x - x \cos x}{x^5}, \quad  R(x) := \frac{1}{\sin^2 x}, \quad T(x) := \frac{x^2 \sin x}{x - \sin x}. $$ As $ R(x) $ is clearly positive and strictly decreasing on $ (0, \pi/2), $ we conclude by Lemmas \ref{lem5} and \ref{lem6} that $ F(x) $ is strictly increasing on $ (0, \pi/2) $ and the required inequality follows.  \hfill $ \qed $ 
 \end{remark}

{\it Proof of Theorem \ref{thm2}.} For the left inequality of \eqref{eqn1.4}, let us set $$ g(x) := \frac{\frac{\sin x}{x} - \frac{2 + \cos x}{3}}{\sin x - x \cos x} = \frac{3 \sin x - x \cos x -2x}{3x \sin x - 3x^2 \cos x}, \quad x \in (0, \pi/2). $$ We now introduce the following intermediary functions: 
$$ g_1(x) := 3 \sin x - x \cos x - 2x, \ h_1(x) := 3x \sin x - 3x^2 \cos x, $$ $$ g_2(x) := x \sin x + 2\cos x - 2, \quad  h_2(x) :=  3\sin x - 3x \cos x + 3x^2 \sin x, $$ and $$ g_3(x) := x\cos x - \sin x,  \quad  h_3(x) := 9x \sin x + 3x^2 \cos x . $$ Then
$$ g_i(0+) = h_i(0+) = 0 \ (i = 1, 2, 3), \quad \frac{g_i'(x)}{h_i'(x)} = \frac{g_{i+1}(x)}{h_{i+1}(x)} \ (i = 1, 2) $$ and
$$ \frac{g_3'(x)}{h_3'(x)} = \frac{-x \sin x}{9 \sin x + 15x \cos x - 3x^2 \sin x} =- \frac{x}{3 v(x)} $$ with $$ v(x) :=  \frac{5x}{\tan x} - x^2 + 3. $$ The function $ v(x) $ is positive on $ (0, \pi/2) $ because of the obvious relation $\tan x/x > 5 / (x^2 - 3) $ and $ v(x) $ is also strictly decreasing on $ (0, \pi/2). $ Applying Lemma \ref{lem1} repeatedly we get $ g(x) $ strictly decreasing on $ (0, \pi/2). $ Hence, for $x < \pi/2$, 
$$ g(x) > g\left(\frac{\pi}{2}-\right)  =  \frac{2}{\pi} - \frac{2}{3}, $$ giving the left inequality of \eqref{eqn1.4}.

For the right inequality of \eqref{eqn1.4}, consider
\begin{align*}
G(x) := \frac{\frac{\sin x}{x} - \frac{2 + \cos x}{3}}{(\sin x - x \cos x)^2} = \frac{1}{3} \cdot  \left[P(x) \right]^2 \cdot Q(x) \cdot R(x),
\end{align*}
where $$ P(x) := \frac{x^2 \sin x}{\sin x - x \cos x}, \quad  Q(x) := \frac{3 \sin x - 2x - x \cos x}{x^5}, \quad  R(x) := \frac{1}{\sin^2 x}. $$ By Lemmas \ref{lem2} and \ref{lem5}, the function $ \left[P(x) \right]^2 \cdot Q(x)  $ is negative increasing and $ R(x) $ is positive decreasing. Therefore, $ G(x) $ is strictly increasing on $ (0, \pi/2) $ and, for $x < \pi/2$, we have $$ G(x) < G\left(\frac{\pi}{2}-\right)  =\frac{2}{\pi} - \frac{2}{3}. $$ This completes the proof. \hfill $ \qed $

\section{Complements}\label{sec4}
This section discusses the importance of our results, by comparing the obtained bounds through analytical and graphical approaches, and also providing alternative proofs. 

\subsection{Comparison of the obtained bound}
The bounds obtained in Theorems \ref{thm1} and \ref{thm2} are comparable. Indeed, we claim that 
\begin{itemize}
\item the lower bound of Theorem 1 is better to the one of Theorem 2,
\item the upper bound of Theorem 2 is better to the one of Theorem 1,
\end{itemize}
which is strictly equivalent to say that $\Phi_1(x) <\Psi_1(x) $ for $x\in (0,\pi/2)$.  After basics developments, we arrive at that is 
\begin{align*}
& \Phi_1(x) <\Psi_1(x) \ \Leftrightarrow \  \frac{2}{\pi}+ \left(1-\frac{2}{\pi}\right)  \cos(x)  < \frac{\sin x}{x},
\end{align*}
which is an inequality proved by \cite{bagul2}. Finally, a  best of Theorems \ref{thm1} and \ref{thm2} is the following inequality: 

\begin{align*}\label{eqn1.3}
\frac{2 + \cos x}{3} - \left(\frac{2}{3} - \frac{2}{\pi}\right)\Phi_1(x) < \frac{\sin x}{x} < \frac{2 + \cos x}{3} - \left(\frac{2}{3} - \frac{2}{\pi}\right)\Psi_2(x). 
\end{align*}

\subsection{Graphical analysis}
\begin{figure}[H]
\begin{center}
\includegraphics[width=12cm,height=7cm]{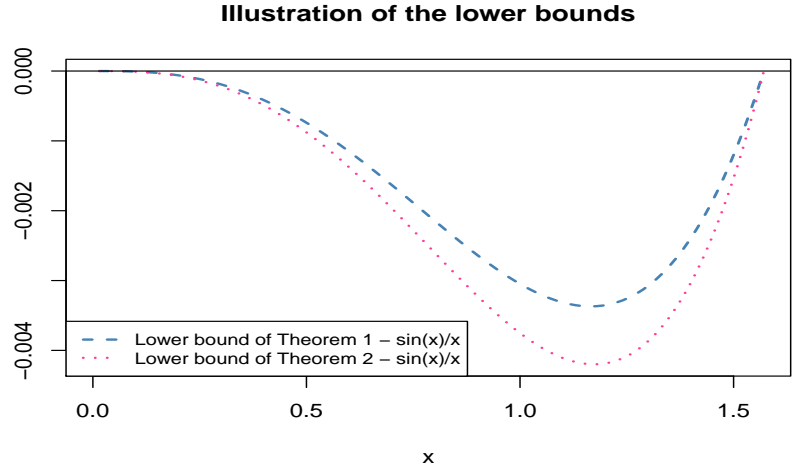}
\caption{Plots of ``lower bounds of Theorems 1 and 2 $-\sin x/x$''}\label{fig1}
\end{center}
\end{figure}

We now provide a graphical analysis of the bounds in Theorems \ref{thm1} and \ref{thm2}, by distinguishing lower bounds and upper bounds. 
Figure \ref{fig1} presents the curves of the lower bounds of the two theorems, both minus  $\sin x/x$ for visual comfort. 

Two immediate remarks come from Figure \ref{fig1}. First, the obtained lower bounds are very sharp; the worst magnitude of the worst of the two curves being a remarkable $\approx 0.004$. As a second remark, it is clear that the lower bound of Theorem \ref{thm1} is uniformly the best, as proved previously. 

\begin{figure}[H]
\begin{center}
\includegraphics[width=12cm,height=7cm]{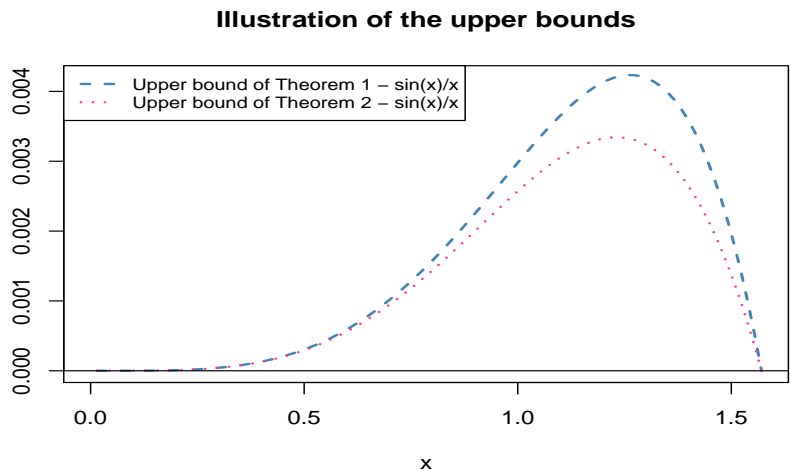}
\caption{Plots of ``upper bounds of Theorems 1 and 2 $-\sin x/x$''}\label{fig2}\end{center}
\end{figure}

Figure \ref{fig2} displays the curves of the upper bounds of the two theorems, both minus  $\sin x/x$. 
 Figure \ref{fig2} illustrates the sharpness of  the obtained upper bounds; the worst magnitude of the worst of the two curves being $\approx 0.004$. It also shows that the upper bound of Theorem \ref{thm2} is uniformly the best, as discussed previously.

\subsection{Alternative proofs of main results} 
It is worth mentioning that alternative proofs of Theorem \ref{thm1} and \ref{thm2} can be developed. Such an alternative is proposed below, based on the methodology developed by \cite{B_Malesevic_M_Makragic_JMI_2016} and \cite{B_D_Banjac_2019}. First, let us notice that the functions
$$
\Phi_1(x) = \dfrac{x - \sin x}{\pi/2 - 1},\quad 
\Phi_2(x) = \left(\dfrac{x - \sin x}{\pi/2 - 1}\right)^{\!2}
$$
and
$$
\Psi_1(x) = \sin x - x \cos x,\quad 
\Psi_2(x) = \left(\sin x - x \cos x\right)^{2}
$$
involved in Theorems \ref{thm1} and \ref{thm2} are Mixed Trigonometric Polynomial (MTP) functions in regards to \cite{B_Malesevic_M_Makragic_JMI_2016} and  \cite{B_D_Banjac_2019}. Also, our main results can be formulated as the following inequality: 
$$
\dfrac{2 + \cos x}{3}  - a \cdot \varphi(x)
<
\frac{\sin x}{x}\,x
<
\dfrac{2 + \cos x}{3}  - a \cdot \psi(x),
$$
where $a =2/3  - 2/\pi$, $\varphi \in \left\{\Phi_1,\Psi_1\right\}$
and $\psi \in \left\{\Phi_2,\Psi_2\right\}$ for $x \in (0,\pi/2)$. Therefore, the previous inequality can be written under the form of an equivalent double sided inequality as
$$
-\varphi(x) < \dfrac{1}{a}  \left( \frac{\sin x}{x}\,x - \dfrac{2 + \cos x}{3} \right) < -\psi(x),
$$
for $x \in (0,\pi/2)$. Since $\varphi(x)$  and $\psi(x)$ are MTP functions, then the problem of proving the previous double sided inequality can be considered by using methods presented in \cite{B_Malesevic_M_Makragic_JMI_2016} and \cite{B_D_Banjac_2019}.
Here, is shown how the considered inequalities in Theorem \ref{thm1} and \ref{thm2} can be reduced to corresponding MTP inequalities.
\begin{itemize}
\item{Case 1: $\varphi(x) = \Phi_1(x) = (x - \sin x)/(\pi/2 - 1).$}
In this case, for $x \in (0,\pi/2)$, we have 
\begin{align*}
& \dfrac{1}{a}\left(\frac{\sin x}{x} x - \dfrac{2+\cos x}{3}\right) + \varphi(x) > 0            \\ 
& \Leftrightarrow  
\dfrac{3 \pi}{2 \pi-6}\left(\frac{\sin x}{x} x - \dfrac{2 + \cos x}{3}\right)
+
\dfrac{x - \sin x}{\pi/2 - 1}
>
0. 
\end{align*}
This is equivalent to prove that $f_1(x)>0$, where 
  \begin{align*}
      f_1(x)& 
:=
\left(\left(-4 \pi+12 \right)x+3 {\pi}^{2}-6 \pi\right)\sin x  \\
&                    
 - 
\pi  \left( \pi-2 \right)   x \cos x
+
x \left(  \left( 4 \pi-12 \right) x-2 {\pi}^{2}+4 \pi \right).
\end{align*}
Or $f_1(x)>0$ is a MTP inequality which can be proved using methods from \cite{B_Malesevic_M_Makragic_JMI_2016} and \cite{B_D_Banjac_2019}.
For example, let us state proof in short form according to \cite{B_D_Banjac_2019}. Let us set 
  \begin{align*}
P_1(x)  := (-4 \pi + 12) x + 3 \pi^2 - 6 \pi,                                       \quad P_2(x)  := - \pi (\pi - 2) x
\end{align*}
and 
  \begin{align*}
P_3(x)  := x \left( (4\pi - 12) x - 2\pi^2 + 4\pi \right).
\end{align*}
Let also $T^{f,0}_{n}(x)$ be {\sc Taylor}'s polynomial function $f(x)$ of $n^{th}$ degree at the point $a=0$.  
Then, it is possible to check that from one side it is true that

\begin{align*}
f_1(x)
& > 
\underbrace{\mathop{P_1(x)}}\limits_{(> 0)}   T^{\sin,0}_7(x)
+
\underbrace{\mathop{P_2(x)}}\limits_{(< 0)}   T^{\cos,0}_4(x) + P_3(x)                          \\
& = 
\left(\frac{\pi}{1260} - \frac{1}{480}\right)  x^8
+
\left(-\frac{\pi^2}{1680} + \frac{1}{840}\right)  x^7                                  \\ 
& + 
\left(-\frac{\pi}{30} + \frac{1}{100}\right)  x^6
+
\left(-\frac{\pi^2}{60} + \frac{\pi}{30}\right)  x^5
+
\left(\frac{2\pi}{3} - 2\right)  x^4                                                        \\ 
& = 
1.12375\ldots \cdot 10^{-4}   x^8
-
2.13477\ldots \cdot 10^{-3}   x^7                                                         \\
& - 
4.71975\ldots \cdot 10^{-3}   x^6
-
5.97736\ldots \cdot 10^{-2}   x^5
+
9.43951\ldots \cdot 10^{-2}   x^4   >                                                          \\ 
& > 
1.12 \cdot 10^{-4}   x^8
-
2.14 \cdot 10^{-3}   x^7                                                           \\ 
& - 
4.72 \cdot 10^{-3}   x^6
-
5.98 \cdot 10^{-2}   x^5
+
9.43 \cdot 10^{-2}   x^4                                                                         \\ 
& = 
\frac{7}{62500}   x^8
-
\frac{107}{50000}   x^7
-
\frac{59}{12500}   x^6
-
\frac{299}{5000}   x^5
+
\frac{943}{10000}   x^4
  >  
0
\end{align*}
for $x \in (0,1.35)$, and that, from other side, it is true that
\begin{align*}
f_1  \left(\mbox{\small $\displaystyle\frac{\pi}{2}$} - x\right)
& > 
\underbrace{\mathop{P_1\left(\mbox{\small $\displaystyle\frac{\pi}{2}$} - x\right)}}\limits_{(> 0)}
  T^{\cos,0}_2(x)
+
\underbrace{\mathop{P_2\left(\mbox{\small $\displaystyle\frac{\pi}{2}$} - x\right)}}\limits_{(< 0)}
  T^{\sin,0}_1(x)
+
P_3\left(\mbox{\small $\displaystyle\frac{\pi}{2}$} - x\right)                              \\ 
& = 
\left(-2 \pi + 6\right)  x^3
+
\left(\pi^2/2 + 2\pi - 12\right)  x^2           + 
\left(-\pi^3/2 - \pi^2 + 12\pi - 12\right) x                                          \\
& = 
-
2.83185\ldots \cdot 10^{-1}   x^3
-
7.82012\ldots \cdot 10^{-1}   x^2
+
3.26369\ldots \cdot 10^{-1}   x                                                                  \\ 
& > 
-
2.84 \cdot 10^{-1}   x^3
-
7.83 \cdot 10^{-1}   x^2
+
3.26 \cdot 10^{-1}   x                                                                           \\ 
& = 
-\frac{71}{250}   x^3 - \frac{783}{1000}   x^2 + \frac{163}{500}   x>
0
\end{align*}
for $x \in (0,0.36)$.

Thus, the MTP inequality $f_1(x) > 0$ for  $x \in (0,\pi/2)$ is wholesomely proved. The other cases are mentioned below, with less details.  

\item{Case 2: $\psi(x) = \Phi_2(x) = \left( (x - \sin x)/(\pi/2 - 1)\right)^{2}.$}
In this case, for $x \in (0,\pi/2)$, we have 
\begin{align*}
& \dfrac{1}{a}\left(\frac{\sin x}{x} x - \dfrac{2+\cos x}{3}\right) + \psi(x) < 0            \\ 
& \Leftrightarrow 
\dfrac{3 \pi}{2 \pi-6}\left(\frac{\sin x}{x} x - \dfrac{2 + \cos x}{3}\right)
+
\left(\dfrac{x - \sin x}{\pi/2 - 1}\right)^{2}
<
0                                                                                                 \\
& \Leftrightarrow 
\mbox{$f_2(x)
:=
-
8 x \left(\pi-3\right) \cos^2 x
-
\pi x \left(\pi-2\right)^{2}\cos x
$}                                                                                               \\ 
&                      +
\left(
\left(-16 \pi+48 \right)x^{2}
+3 \pi\left(\pi - 2\right)^2 \right) \sin  x                                                                                               \\ 
&                     
 -  
2 {\pi}^{3}x+8 \pi x^{3}+8 {\pi}^{2}x-24 {x}^{3}-24 x
<
0.
\end{align*}
The MTP inequality $
f_2(x) < 0 $
over $(0, \pi/2)$ can be proven using methods from \cite{B_Malesevic_M_Makragic_JMI_2016} and \cite{B_D_Banjac_2019}.

\item{Case 3: $\varphi(x) = \Psi_1(x) = \sin x - x \cos x.$}
In this case, for $x \in (0,\pi/2)$, it holds that
\begin{align*}
& \dfrac{1}{a}\left(\frac{\sin x}{x} x - \dfrac{2+\cos x}{3}\right) + \varphi(x) > 0            \\ 
& \Leftrightarrow 
\dfrac{3 \pi}{2 \pi-6}\left(\frac{\sin x}{x} x - \dfrac{2 + \cos x}{3}\right)
+
\sin x - x \cos x
>
0                                                                                                 \\ 
& \Leftrightarrow  
 f_3(x)
:=
\left(2\left(\pi-3\right)+3\pi\right)\sin x      +  
x \left(-2\left(\pi-3\right)x-\pi\right) \cos x - 2 \pi x
>
0. 
\end{align*}
Hence, the MTP inequality $
f_3(x) > 0$
over $(0, \pi/2)$ can be proved via the methods from
\cite{B_Malesevic_M_Makragic_JMI_2016} and \cite{B_D_Banjac_2019}.

\item{Case 4: $\psi(x) = \Psi_2(x) = \left(\sin x - x \cos x\right)^2.$} For $x \in (0,\pi/2)$, we have
\begin{align*}
&                    
\dfrac{1}{a}\left(\frac{\sin x}{x} x - \dfrac{2+\cos x}{3}\right) + \psi(x) < 0            \\ 
& \Leftrightarrow 
\dfrac{3 \pi}{2 \pi-6}\left(\frac{\sin x}{x} x - \dfrac{2 + \cos x}{3}\right)
+
\left(\sin x - x \cos x\right)^2
<
0                                                                                                 \\ 
& \Leftrightarrow    f_4(x)
:=
-
4 x^2 \left(\pi-3\right)\sin x \cos x
+
3 \pi \sin x \\
&            +  
2 x (x^2-1)(\pi-3) \cos^2 x - \pi x \cos x - 6 x
<
0. 
\end{align*}
And, the 
MTP inequality
$f_4(x) < 0$
over $(0, \pi/2)$ can be proven using methods from \cite{B_Malesevic_M_Makragic_JMI_2016} and \cite{B_D_Banjac_2019}.

\end{itemize}
This ends this possible alternative proof.

\break

\end{document}